\documentclass{article}
\usepackage{longtable}
\usepackage{booktabs}
\usepackage{multirow}
\usepackage{scalefnt}
\usepackage{amsmath,amsthm,amsfonts,amssymb}
\usepackage{bm}
\usepackage{here}
\usepackage{caption}
\usepackage{secdot}

\sectiondot{subsection}
\usepackage[margin=10truemm]{geometry}
\usepackage{url}
\usepackage{cases}
\usepackage{empheq}
\numberwithin{equation}{section}
\title{\huge The problem of finding three numbers such that the sum or difference of any two of them is a square number}
\author{Seiji Tomita }
\date{}
\begin{document}
\maketitle
\begin{abstract}
Euler explored the problem of finding three numbers such that the sum or difference of any two of them is a perfect square. He discovered a parametric solution represented by polynomials of degree 18 and identified the smallest of these solutions. Additionally, we derive two parametric solutions using the elliptic curve method and find all solutions within the range of $10^8$.
\end{abstract}
\vskip3\baselineskip

\section{ Introduction}
\vskip3\baselineskip
The Six Squares Problem, also known as Mengoli's Six-Square Problem,
 is a mathematical problem that explores the relationship between squares
 and integers. It involves finding integers such that their differences are
 squares and the differences of their squares are also squares.
The problem asks for three integers, $x,\ y,$ and $z,$ 
where the differences $(x-y,\ x-z,\ y-z)$ are all perfect squares,
 and the differences of their squares $(x^2-y^2,\ x^2-z^2,\ y^2-z^2)$ are also perfect squares.
  
Ozanam\cite{c} found the solution $(x, \ y, \ z)=(2399057 ,\  2288168 ,\  1873432)$ in $1691.$
A little later, Euler\cite{a} found the smallest solution $(x,y,z)=(434657 ,  420968 ,  150568)$ in $1770.$
Mengoli\cite{c} obtained the solution $(x,\ y,\ z)=( 40606322,\ 29316722,\ 26633678  ).$

In this paper, we show that Euler's method and two new approaches using elliptic curves produce two parametric solutions represented by polynomials of degree 20.
These two new methods generate an infinite number of parametric solutions.

Moreover, we have re-verified that Euler's solution \((x, y, z) = (434657, 420968, 150568)\) is indeed the smallest solution through a brute-force search. Additionally, since the numerical solutions derived from the polynomials are large, we employed brute force to obtain smaller solutions.

Given that his paper is written in Latin, it is essential to translate it into English 
. This will not only make the work accessible to wider researchers but also highlight its significance and contributions.

\vskip3\baselineskip

\section{Euler's method}
\vskip\baselineskip
Finding integers $x,y,$ and $z$ such that $x \pm y, x \pm z, y \pm z$ are square numbers.
\vskip\baselineskip

Euler obtained the simple parametric solution and the smallest positive integer solution.
First, he considers the case $x=p^2+q^2=r^2+s^2.$
Let $x=p^2+q^2$ and $y=2pq.$ Then $(x+y,\ x-y)=( (p+q)^2,\ (p-q)^2 ).$ 
Similarly, let $x=r^2+s^2$ and  $z=2rs.$ We get $(x+z,\ x-z)=( (r+s)^2,\ (r-s)^2 ).$
\ To satisfy that $p^2+q^2=r^2+s^2,$ he used the well--nown parametric solution as follows.

$$(p,\ q, \ r, \ s)=( ac + bd,\  ad - bc, \ ad + bc,\  ac - bd).$$

Thus, we have $$(y+z,\ y-z)=(2pq+2rs,\ 2pq-2rs)=(4cd(a-b)(a+b),\ 4ab(d-c)(d+c)).$$

Since  $y+z$ and $y-z$ are both perfect squares, we can express the equation $$y^2-z^2=16acbd(a^2-b^2)(d^2-c^2)$$

 as a perfect square as follows.
\begin{equation}
acbd(a^2-b^2)(d^2-c^2)=v^2.
\end{equation}

Next, we find the rational soution for $(b,\ c)$ such that
$$cd(d^2-c^2)=n^2ab(a^2-b^2)$$

 where $d=a$ and $a=c-b.$

Then, we have

$$\frac{b}{c} = \frac{2+n^2}{1+2n^2}.$$

Hence, we take $$(b,\ c)=(2+n^2,\ 1+2n^2).$$

Since $ab(d^2-c^2)$ must be a square number, we can reduce it to $$ab(d^2-c^2)=3(1-n^2)(2+n^2)^2n^2.$$

Then, $3(1-n^2)$ must be a square number.

We solve  $$3(1-n^2)=\frac{f^2}{g^2}(n+1)^2$$

 for $n$ where $f$ and $g$ are any inetgers, we get

$$n = -\frac{-3g^2+f^2}{3g^2+f^2}.$$

Thus, we obtain $(a,\ b,\ c),$ 

\begin{align*}
a &=-4f^2g^2, \\
b &=9g^4+2g^2f^2+f^4, \\
c &=9g^4-2g^2f^2+f^4.
\end{align*}

Finally, he obtained the parametric solution of 16 degree.

\begin{align*}
x &= f^{16}+60g^4f^{12}+1318g^8f^8+4860g^{12}f^4+6561g^{16}, \\
y &= 16g^2f^{14}+112g^6f^{10}+1008g^{10}f^6+11664g^{14}f^2, \\
z &= 32g^4f^{12}+960g^8f^8+2592g^{12}f^4.
\end{align*}

He gave the example for $(f,\ g)=(2,\ 1).$

$$(x,\ y,\ z)=(733025,\ 488000,\ 418304).$$

$$(x+y,\ y+z,\ x+z,\ x-y,\ y-z,\ x-z)=(1105^2,\ 952^2,\ 1073^2,\ 495^2,\ 264^2,\ 561^2).$$

\begin{table}[hbtp]
\centering
  \caption{Numerical solutions with $(f,g)<5 $}
  \label{tab:table1}
  \begin{tabular}{rrrrrrrrrrrrrr}
    \hline
   $f$ & $g$  & $x$ & $y$  & $z$ &  $\sqrt{x+y}$ &  $\sqrt{x-y}$ &  $\sqrt{x+z}$ &  $\sqrt{x-z}$ &  $\sqrt{y+z}$ &  $\sqrt{y-z}$ \\
    \hline \hline
  1 &  2  &  450226625 &  192142400 &  10863104  &   25345 &  16065 &  21473 &  20961 &  14248 &  13464 \\
  1 &  3  &  1113363250 &  218156850 &  5405454  &  36490 &  29920 &  33448 &  33286 &  14952 &  14586 \\
  1 &  4  &  28260904090625 &  3132088582400 &  43549466624  &  5602945 &  5012865 &  5320193 &  5312001 &  1782032 &  1757424 \\
  2 &  1  &  733025 &  488000 &  418304  &  1105 &  495 &  1073 &  561 &  952 &  264 \\
  2 &  3  &  325988094625 &  227049537600 &  23662969344  &  743665 &  314545 &  591313 &  549841 &  500712 &  450984 \\
  3 &  4  &  5389227425 &  4416723200 &  600449024  &  99025 &  31185 &  77393 &  69201 &  70832 &  61776 \\
  4 &  3  &  1596175283425 &  1260685267200 &  808907710464  &  1690225 &  579215 &  1550833 &  887281 &  1438608 &  672144 \\
    \hline
  \end{tabular}
\end{table}
\vskip2\baselineskip
Since Euler's parametric solution gives only a partial solution  of equation $(2.1),$ we searched the positive integer solutions of equation $(2.1)$ where $(a,b,c,d)<300$ and $(x,y,z)<10^8. $

\begin{table}[hbtp]
\centering
  \caption{The small solutions with $(x,y,z)<10^8 $}
  \label{table:table2}
  \begin{tabular}{rrrrrrrrrrrrrr}
    \hline
   $a$ & $b$ & $c$ & $d$ & $x$ & $y$  & $z$ &  $\sqrt{x+y}$ &  $\sqrt{x-y}$ &  $\sqrt{x+z}$ &  $\sqrt{x-z}$ &  $\sqrt{y+z}$ &  $\sqrt{y-z}$ \\
    \hline \hline
 26 &  10 &  49 &  81 & 434657 &  420968 &  150568 & 925 &  117 &  765 &  533 &  756 &  520 \\
 17 &  1 &  32 &  49 & 993250 &  949986 &  856350 & 1394 &  208 &  1360 &  370 &  1344 &  306 \\
 33 &  1 &  17 &  49 & 733025 &  488000 &  418304 & 1105 &  495 &  1073 &  561 &  952 &  264 \\
 200 &  50 &  240 &  289 & 2399057 &  2288168 &  1873432 & 2165 &  333 &  2067 &  725 &  2040 &  644 \\
 9 &  4 &  117 &  125 & 2843458 &  2040642 &  1761858 & 2210 &  896 &  2146 &  1040 &  1950 &  528 \\
 8 &  7 &  121 &  135 & 3713858 &  891458 &  88642 & 2146 &  1680 &  1950 &  1904 &  990 &  896 \\
 49 &  39 &  49 &  55 & 5320193 &  1782032 &  589568 & 2665 &  1881 &  2431 &  2175 &  1540 &  1092 \\
 28 &  17 &  55 &  64 & 7640833 &  4504392 &  2465208 & 3485 &  1771 &  3179 &  2275 &  2640 &  1428 \\
 39 &  1 &  95 &  121 & 9004913 &  8845712 &  8626688 & 4225 &  399 &  4199 &  615 &  4180 &  468 \\
34 &  14 &  121 &  240 & 24417458 &  24163442 &  3714958 &  6970 &  504 &  5304 &  4550 &  5280 &  4522 \\
41  &   1  &   64  &   105   &   25433522  &   23147378  &   22011022   &   6970  &   1512  &   6888  &   1850  &   6720  &   1066 \\
17  &   4  &   256  &   273   &   42719825  &   39381896  &   36935800   &   9061  &   1827  &   8925  &   2405  &   8736  &   1564 \\
28  &   1  &   87  &   256   &   57387425  &   38124104  &   31631800   &   9773  &   4389  &   9435  &   5075  &   8352  &   2548 \\
81  &   49  &   117  &   125   &   65678017  &   34261992  &   26578008   &   9997  &   5605  &   9605  &   6253  &   7800  &   2772 \\
    \hline
  \end{tabular}
 \vspace{20mm}
\end{table}

\newpage
\section{Elliptic curves method I}

\vskip\baselineskip
We tackle the problem using a different approach to Euler's method and show how infinite parametric solutions can be derived.
\vskip\baselineskip

First, we set $$(x,\ y,\ z)=(m^2n^2+1,\ m^2+n^2,\ 2mn).$$

Then,  $(x  \pm y,\ x  \pm z,\ y  \pm z)$ is written as
$$(x+y,\ x+z,\ y+z)=((n^2+1)(1+m^2),\ (mn+1)^2,\ (m+n)^2),$$

$$(x-y,\ x-z,\ y-z)=((n-1)(n+1)(m-1)(m+1),\ (mn-1)^2,\ (m-n)^2).$$

Thus, $x \pm y$ must be square numbers.

First, we parametrise $x+y=m^2n^2+1+m^2+n^2=u^2$ using a known solution $(m,\ u)=(n,\ n^2+1).$

We get $m=-\dfrac{n^3+n+nk^2-2k-2n^2k}{n^2+1-k^2},$ and substitute it into $x-y=v^2.$

Let $V = v(-1-n^2+k^2)$ and $U=k.$ Then we have a quartic equation 

\begin{equation}
Q: V^2 = (n^2-1)^2U^4+(4n-4n^5)U^3+(-6-6n^2+6n^4+6n^6)U^2+(4n+4n^3-4n^5-4n^7)U+(n^4-1)^2.
\end{equation}

$(3.1)$ is in three variables $n, U,$ and $V$ , and represents an algebraic
surface over the field $\mathbb{Q}$ , we  regard it as a quartic model of over the function field $\mathbb{Q}(n)$.

$(3.1)$ is reduced to the following cubic equation.
(3.2)  represents an algebraic surface over $\mathbb{Q}$ , we regard it as the Weierstrass model of an elliptic curve over the function
field $\mathbb{Q}(n)$.

\begin{equation}
E: Y^2+a1XY+a3Y=X^3+a2X^2+a4X+a6
\end{equation}

where
\begin{align*}
a1 &=-4n^3-4n,\\
a2 &=2n^6-2n^4-10n^2-6,\\
a3 &=16n^5-8n-8n^9,\\
a4 &=8n^2-16n^6+8n^{10}-4+4n^4+4n^8-4n^{12},\\
a6 &=-8n^{18}+24n^{16}+32n^{14}-96n^{12}-48n^{10}+144n^8+32n^6-96n^4-8n^2+24.\\
\end{align*}
The transformation from $(3.2)$ to $(3.1)$ is given
\begin{align*}
U = &(4n^{10}-4n^8-24n^6+2n^4X-8n^4+20n^2-2X+12)/Y,\\
V = &(24-76X+80n^{14}+8n^{22}-24n^{20}-240n^{12}-168n^4X-X^3+18X^2\\
    &+320n^8X-80n^{10}+120n^{16}+240n^8-8n^2-120n^4+40n^6-40n^{18}+32nY\\
    &+12n^{16}X-232Xn^2+184Xn^6+72Xn^{10}-24Xn^{14}-88Xn^{12}+30n^2X^2\\
    &-36n^6X^2-12n^4X^2-6n^8X^2+X^3n^4+6X^2n^{10}+64n^3Y-32n^9Y-64n^7Y)/Y^2.
\end{align*}

Since $E$ has a known point $P=(-2n^6+2n^4+10n^2+6,\ 32n^5+64n^3+32n),$

we get the rational point $2P$ using the group law as follows.

$$2P(X,\ Y)=(-(n^8+2n^6-4n^4-2n^2-1)/n^2,\ (n^{12}-6n^{10}+n^8+12n^6+15n^4+2n^2-1)/n^3)$$
\vskip\baselineskip

As the point $2P$ does not have integer coordinates, the point $P$ has infinite order by the Nagell-Lutz theorem\cite{d}.
Since there are infinitely many rational points on the elliptic curve $E$ and these can be obtained by
the group law.
Thus, we can obtain infinitely many rational solutions of $(3.1).$

For instance, we will get the rational point $2Q(U, \ V)$.
We can pull back to  the  corresponding  point $2Q(U, \ V)$ of the quartic equation below.
\begin{align*}
U &= \frac{2n(n-1)(n+1)(n^2+1)}{n^4-2n^2-1},\\
V &= \frac{(n-1)(n+1)(n^2+1)(n^4-2n^3-2n+1)(n^4+2n^3+2n+1)}{(n^4-2n^2-1)^2}.
\end{align*}

Hence, substitute $k = \dfrac{2n(n-1)(n+1)(n^2+1)}{n^4-2n^2-1}$ into $(x,\ y,\ z)$, we get
\begin{align*}
x &= (n^4+1)(n^{16}+20n^{12}-26n^8+20n^4+1),\\
y &= 2n^2(5n^8-2n^4+1)(n^8-2n^4+5),\\
z &= 2n^2(3n^8-6n^4-1)(n^8+6n^4-3).
\end{align*}
\vskip\baselineskip

Indeed, $x \pm y, x \pm z, y \pm z$ are square numbers.
\begin{align*}
x+y &= (n^2+1)^2(n^8+4n^6-6n^4+4n^2+1)^2,\\
x+z &= (n^2-1)^2(n^4-2n^3+4n^2-2n+1)^2(n^4+2n^3+4n^2+2n+1)^2,\\
y+z &= 16n^2(n^2-1)^2(n^2+1)^2(n^4+1)^2,\\
x-y &= (n^2-1)^2(n^4-2n^3-2n+1)^2(n^4+2n^3+2n+1)^2,\\
x-z &= (n^2+1)^2(n^4-2n^3+2n+1)^2(n^4+2n^3-2n+1)^2 ,\\
y-z &= 4n^2(n^4-2n^2-1)^2(n^4+2n^2-1)^2.
\end{align*}
\vskip\baselineskip

For instance, we get $(x,\ y,\ z)=(2399057,\ 2288168,\ 1873432)$ for $n=2.$
\vskip2\baselineskip

We searched the rational solutions of equation $(3.1)$ where $n<2000$ and $(x,\ y,\ z)<10^8. $
\begin{table}[hbtp]
\centering
  \caption{The small solutions of equation $(3.1)$with $(x,\ y,\ z)<10^8 $}
  \label{table:table3}
  \begin{tabular}{rrrrrrrrrrr}
    \hline
   $n$ & $m$ &  $x$ & $y$  & $z$ &  $\sqrt{x+y}$ &  $\sqrt{x-y}$ &  $\sqrt{x+z}$ &  $\sqrt{x-z}$ &  $\sqrt{y+z}$ &  $\sqrt{y-z}$ \\
    \hline \hline
2/11 &  58/59 & 434657 &  420968 &  150568 &  925 &  117 &  765 &  533 &  756 &  520 \\
8/19 &  32/43 & 733025 &  488000 &  418304 &  1105 &  495 &  1073 &  561 &  952 &  264 \\
3/5 &  165/173 & 993250 &  949986 &  856350 &  1394 &  208 &  1360 &  370 &  1344 &  306 \\
2 &  698/671 & 2399057 &  2288168 &  1873432 &  2165 &  333 &  2067 &  725 &  2040 &  644 \\
7/9 &  79/177 & 2843458 &  2040642 &  1761858 &  2210 &  896 &  2146 &  1040 &  1950 &  528 \\
1/41 &  23/47 & 3713858 &  891458 &  88642 &  2146 &  1680 &  1950 &  1904 &  990 &  896 \\
4/7 &  32/329 & 5320193 &  1782032 &  589568 &  2665 &  1881 &  2431 &  2175 &  1540 &  1092 \\
2/9 &  226/303 & 7640833 &  4504392 &  2465208 &  3485 &  1771 &  3179 &  2275 &  2640 &  1428 \\
28/29 &  64/83 & 9004913 &  8845712 &  8626688 &  4225 &  399 &  4199 &  615 &  4180 &  468 \\
 11/17 &  229/257 & 25433522 &  23147378 &  22011022 & 6970 &  1512 &  6888 &  1850 &  6720 &  1066 \\
 10/11 &  326/515 & 42719825 &  39381896 &  36935800 & 9061 &  1827 &  8925 &  2405 &  8736 &  1564 \\
 5/2 &  1451/1090 & 57387425 &  38124104 &  31631800 & 9773 &  4389 &  9435 &  5075 &  8352 &  2548 \\
 2/3 &  838/2643 & 65678017 &  34261992 &  26578008 & 9997 &  5605 &  9605 &  6253 &  7800 &  2772 \\
 49/43 &  167/29 & 68516498 &  53586002 &  20408402 & 11050 &  3864 &  9430 &  6936 &  8602 &  5760 \\
 7 &  1249/1151 & 77764850 &  66475250 &  20126386 & 12010 &  3360 &  9894 &  7592 &  9306 &  6808 \\
 8/15 &  400/561 & 81052225 &  56142144 &  53856000 & 11713 &  4991 &  11615 &  5215 &  10488 &  1512 \\
    \hline
  \end{tabular}
    \vspace{20mm}
\end{table}
\vskip5\baselineskip

\newpage
\section{Elliptic curves method I\hspace{-1.2pt}I}
\vskip\baselineskip
We tackle the problem differently than we did in the previous section.

Let us solve the following six simultaneous equations in integers:
\begin{empheq}[left=\empheqlbrace]{align}
           x+y&=p^2  \\
           y+z&=q^2  \\
           x+z&=r^2  \\
           x-y&=u^2  \\
           y-z&=v^2  \\
           x-z&=w^2  
\end{empheq}

Using equations $(4.1),(4.2),$ and $(4,3),$ we obtain $(x,\ y,\ z),$

$$(x,\ y,\ z)= \Bigl( \frac{p^2-q^2+r^2}{2},\ \frac{p^2+q^2-r^2}{2},\ \frac{-p^2+q^2+r^2}{2} \Bigr).$$

Substitute $(x, \ y,\ z)$ into $(4.4),(4.5),$ and $(4.6),$ we get

 \begin{align}
r^2-q^2=u^2 &, \\
p^2-r^2=v^2 & ,\\
p^2-q^2=w^2 &.
 \end{align}

From equation (4.7), let $(r,\ q,\ u)$ be defined as $(m^2 + n^2,\ m^2 - n^2,\ 2mn).$

Equations $(4.8)$ and $(4.9)$ yield $v^2-w^2 = q^2-r^2 =-4m^2n^2.$

Solving for $(v,\ w),$ we get
$$(v,\ w)=\Bigl(\frac{n^2-m^2t^2}{t},\ \frac{m^2t^2+n^2}{t} \Bigr)$$

 where $t$ is any rational number.
 
Thus, equation $(4.9)$ reduces to $$p^2 = (\frac{1}{t^2}+1)n^4+(t^2+1)m^4.$$

Hence, we have the quartic equation
 \begin{equation}
Q: V^2 = (t^2+1)U^4+(t^2+1)t^2
 \end{equation}

where $$(U,\ V)=(\dfrac{n}{m},\ \dfrac{tp}{m^2}).$$
As in the previous section, we  regard it as a quartic model of over the function field $\mathbb{Q}(t)$.

Since the quartic equation $(4.10)$ has a rational point $Q(U,\ V)=(t,\ t(t^2+1)),$ then this quartic equation is birationally equivalent to an elliptic curve below.
$(4.10)$ is reduced to the following cubic equation.
$(4.11)$  represents an algebraic surface over $\mathbb{Q}$ , we regard it as the Weierstrass model of an elliptic curve over the function
field $\mathbb{Q}(t)$.

\begin{align}
E: Y^2+4t^2YX+(8t^6+16t^4+8t^2)Y &= X^3+(2t^4+6t^2)X^2 \\
                                        &+(-4t^8-12t^6-12t^4-4t^2)X  \notag \\
                                        &-8t^{12}-48t^{10}-96t^8-80t^6-24t^4. \notag
\end{align}

Since $E$ has a known point $P=(-2t^4-6t^2,\ 8t^4-8t^2)$, we can obtain the rational point $2P$ using the group law as follows.

$$2P(X,\ Y)= \Bigl( \frac{(t^4-6t^2-3)t}{3t^4+6t^2-1},\ \frac{(t^2+1)(t^8+28t^6+6t^4+28t^2+1)t}{(3t^4+6t^2-1)^2} \Bigr)$$

We can pull back to  the  corresponding  point $2Q(U, \ V)$ of the quartic equation below.
The transformation from $(4.11)$ to $(4.10)$ is given

\begin{align*}
U &= (4t^7+16t^5+2t^3X+12t^3+tY+2tX)/Y,\\
V &= (4t^3X+104t^7+18t^3X^2+tX^3+t^3X^3+12t^{11}X+72t^9X+144t^7X+6t^7X^2 \\
    &  -8t^7Y+24t^5X^2+88t^5X+176t^9+144t^{11}+8Yt^3+56t^{13}+8t^{15}+24t^5)/(Y^2).
\end{align*}

Hence, substitute $k = $ into $(p,\ q,\ r)$, we get

\begin{align}
p &= (t^2+1)(t^8+28t^6+6t^4+28t^2+1),\\
q &= t^{10}-21t^8-6t^6+6t^4+21t^2-1,\\
r &= t^{10}-3t^8+66t^6+66t^4-3t^2+1.
\end{align}

\begin{align*}
x &= 2(t^4+6t^2+1)(t^{16}+88t^{14}+92t^{12}+872t^{10}+1990t^8+872t^6+92t^4+88t^2+1),\\
y &= (t^8-8t^7+44t^6-56t^5+102t^4-56t^3+44t^2-8t+1)(t^8+8t^7+44t^6+56t^5+102t^4+56t^3+44t^2+8t+1)(t^2-1)^2,\\
z &= 2(t^8+8t^7-20t^6+56t^5-26t^4+56t^3-20t^2+8t+1)(t^8-8t^7-20t^6-56t^5-26t^4-56t^3-20t^2-8t+1)(t^2-1)^2.
\end{align*}
\vskip\baselineskip
Indeed, $x \pm y,\ x \pm z,\ y \pm z$ are square numbers as follows.
\begin{align*}
x+y &= 4(t^2+1)^2(t^8+28t^6+6t^4+28t^2+1)^2,\\
x+z &= 4(t^2+1)^2(t^4-4t^3+6t^2+4t+1)^2(t^4+4t^3+6t^2-4t+1)^2,\\
y+z &= 4(t-1)^2(t+1)^2(t^4+4t^3-2t^2+4t+1)^2(t^4-4t^3-2t^2-4t+1)^2,\\
x-y &= 16t^2(3t^4+6t^2-1)^2(t^4-6t^2-3)^2,\\
x-z &= 16t^2(5t^4+2t^2+1)^2(t^4+2t^2+5)^2 ,\\
y-z &= 256t^2(t-1)^2(t+1)^2(t^2+1)^2(t^4+6t^2+1)^2.
\end{align*}
\vskip\baselineskip
For instance, we obtain $(x,\ y,\ z)=(434657 ,\  420968 ,\  150568)$ for $t=1/117.$
\vskip\baselineskip
We searched the rational solutions of equation $(4.10)$ where $t<1000$ and $(x,\ y,\ z)<10^8. $
\vskip\baselineskip
$(x,\ y,\ z)=(2843458 ,\  2040642 ,\  1761858)$ and $(3713858 ,\  891458 ,\  88642)$ were not found in the search range.
\newpage
\begin{table}[hbtp]
\centering
  \caption{The small solutions of equation $(4.10)$ with $(x,\ y,\ z)<10^8	 $}
  \label{table:table4}
  \begin{tabular}{rrrrrrrrrr}
    \hline
   $t$ & $x$ & $y$  & $z$ &  $\sqrt{x+y}$ &  $\sqrt{x-y}$ &  $\sqrt{x+z}$ &  $\sqrt{x-z}$ &  $\sqrt{y+z}$ &  $\sqrt{y-z}$ \\
    \hline \hline
  1/117  &  434657 &  420968 &  150568  &  925 &  117 &  765 &  533 &  756 &  520 \\
  11/27  &  733025 &  488000 &  418304  &  1105 &  495 &  1073 &  561 &  952 &  264 \\
  1/3  &  2399057 &  2288168 &  1873432  &  2165 &  333 &  2067 &  725 &  2040 &  644 \\
  3/11  &  5320193 &  1782032 &  589568  &  2665 &  1881 &  2431 &  2175 &  1540 &  1092 \\
  7/11  &  7640833 &  4504392 &  2465208  &  3485 &  1771 &  3179 &  2275 &  2640 &  1428 \\
  57  &  9004913 &  8845712 &  8626688  &  4225 &  399 &  4199 &  615 &  4180 &  468 \\
  1/21  &  42719825 &  39381896 &  36935800  &  9061 &  1827 &  8925 &  2405 &  8736 &  1564 \\
  3/7   &  57387425 &  38124104 &  31631800  &  9773 &  4389 &  9435 &  5075 &  8352 &  2548 \\
  5 &    65678017 &  34261992 &  26578008  &   9997 &  5605 &  9605 &  6253 &  7800 &  2772 \\
  7/23 &    81052225 &  56142144 &  53856000  &   11713 &  4991 &  11615 &  5215 &  10488 &  1512 \\
    \hline
  \end{tabular}
  \vspace{20mm}
\end{table}

\section{Brute-force search}
\vskip\baselineskip

A brute-force search of equation $(4.1) \sim (4.6)$ yields the following result where $(x,\ y,\ z)<10^7.$

We re-verified that $(x,y,z)=(434657 , \ 420968 , \ 150568)$ is the smallest solution.

\begin{table}[hbtp]
\centering
  \caption{A brute-force search with $(x,\ y,\ z)<10^7	 $}
  \label{table:table5}
  \begin{tabular}{rrrrrrrrr}
    \hline
    $x$ & $y$  & $z$ &  $\sqrt{x+y}$ &  $\sqrt{x-y}$ &  $\sqrt{x+z}$ &  $\sqrt{x-z}$ &  $\sqrt{y+z}$ &  $\sqrt{y-z}$ \\
    \hline \hline
 434657 &  420968 &  150568 &  925 &  117 &  765 &  533 &  756 &  520 \\
 733025 &  488000 &  418304 &  1105 &  495 &  1073 &  561 &  952 &  264 \\
 993250 &  949986 &  856350 &  1394 &  208 &  1360 &  370 &  1344 &  306 \\
 2399057 & 2288168 & 1873432 & 2165 & 2040 & 2067 &   333 &  644  & 725 \\
2843458 & 2040642 & 1761858 & 2210 & 1950 & 2146  &  896 &  528 & 1040 \\
3713858  & 891458 & 88642 & 2146 &  990 & 1950 &  1680 &  896    & 1904 \\
5320193  &1782032 & 589568 & 2665 & 1540 & 2431  & 1881 & 1092 & 2175 \\
7640833 &  4504392 &  2465208 & 3485 &  1771 &  3179 &  2275 &  2640 &  1428 \\
9004913 & 8845712 & 8626688 & 4225 & 4180 & 4199 &   399 &  468 &  615 \\
    \hline
  \end{tabular}
\vspace{20mm}
 \end{table}

\section{Connect to the tetrahedron}

A tetrahedron is a polyhedron with four triangular faces, six edges, and four vertices. It's the simplest convex polyhedron and can be seen as a triangular pyramid.
While the problem doesn't directly define the tetrahedron, the solutions to the six-square problem might be used in constructing a tetrahedron or in exploring properties related to its faces, edges, or vertices.  
$(4.1) \sim (4.6)$ yields the following result.
$x^2-y^2=(pu)^2,\ x^2-z^2=(vw)^2,\ y^2-z^2=(qv)^2.$
This gives three Pythagorean triples, $(y,\ pu,\ x), (z,\ vw,\ x)$, and $(z,\ qv,\ y).$  Again, we have a tetrahedron, one with three right triangles.
 The triangle that is not a right triangle has sides of length pu,\ vw, and qv. 
 The three right angles do not share a common vertex.  

Indeed, as shown in Table $3,$ we have
\begin{align*}
x^2-y^2 = 434657^2 -  420968^2 = (925 \cdot 117)^2,\\
x^2-z^2 = 434657^2 -  150568^2 = (765 \cdot 533)^2,\\
y^2-z^2 = 420968^2 -  150568^2 = (756 \cdot 520)^2.
\end{align*}

Incidentally, a different solution of the following system equation  was found in the previous section.
\begin{align*}
\begin{cases}
x^2-y^2=u^2\\
x^2-z^2=v^2\\
y^2-z^2=w^2
\end{cases}
\end{align*}
 
From $(4.7) \sim (4.14),$ taking $x=p,\ y=r,$ and \ $z=q,$ we obtain

\begin{align*}
x &= (t^2+1)(t^8+28t^6+6t^4+28t^2+1),\\
y &= t^{10}-3t^8+66t^6+66t^4-3t^2+1,\\
z &= t^{10}-21t^8-6t^6+6t^4+21t^2-1,\\
u &= 8t(t-1)(t+1)(t^2+1)(t^4+6t^2+12),\\
v &= 2t(t^4+2t^2+5)(5t^4+2t^2+1),\\
w &= 2t(3t^4+6t^2-1)(t^4-6t^2-3).
\end{align*}

Taking $t=3$ and removing the common factor, we get
\begin{align*}
x^2-y^2 &= 4330^2 -  1450^2 = 4080^2,\\
x^2-z^2 &= 4330^2 -  1288^2 = 4134^2,\\
y^2-z^2 &= 1450^2 -  1288^2 = 666^2.
\end{align*}

\vskip3\baselineskip

\section{Concluding Remarks}
\vskip\baselineskip
A brute force search revealed no solutions smaller than $434657,$ confirming once again that Euler's solution is the smallest.

Our method produces solutions $(x,y,z)=(68516498 ,  53586002 ,  20408402), \ (77764850 ,  66475250,  20126386)$ and

 $(81052225 ,  56142144 ,  53856000),$ but not Euler's method.

On the other hand, solution $(x,y,z)=(24417458 ,  24163442 ,  3714958)$ appears in the Euler method, but not in our one.

Thus, finding a parametric solution that generates all integer solutions remains an open problem.

\end{document}